\documentclass[12pt,amsmath]{article}
\usepackage{amsfonts}
\begin{document}
\newsymbol\rtimes 226F
\newsymbol\ltimes 226E
\newcommand{\text}[1]{\mbox{{\rm #1}}}
\newcommand{\Rep}{\text{Rep}}
\newcommand{\Fun}{\text{Fun}}
\newcommand{\Hom}{\text{Hom}}
\newcommand{\End}{\text{End}}
\newcommand{\GL}{\text{GL}}
\newcommand{\Sp}{\text{Sp}}
\newcommand{\Ps}{\text{Ps}}
\newcommand{\ASp}{\text{ASp}}
\newcommand{\APs}{\text{APs}}
\newcommand{\gd}{\delta}
\newcommand{\lan}{\langle}
\newcommand{\ran}{\rangle}
\newcommand{\itms}[1]{\item[[#1]]}
\newcommand{\nin}{\in\!\!\!\!\!/}
\newcommand{\g}{{\bf g}}
\newcommand{\sub}{\subset}
\newcommand{\cntd}{\subseteq}
\newcommand{\go}{\omega}
\newcommand{\Pa}{P_{a^\nu,1}(U)}
\newcommand{\fx}{f(x)}
\newcommand{\fy}{f(y)}
\newcommand{\gD}{\Delta}
\newcommand{\gl}{\lambda}
\newcommand{\gL}{\Lambda}
\newcommand{\half}{\frac{1}{2}}
\newcommand{\sto}[1]{#1^{(1)}}
\newcommand{\stt}[1]{#1^{(2)}}
\newcommand{\Z}{\hbox{\sf Z\kern-0.720em\hbox{ Z}}}
\newcommand{\singcolb}[2]{\left(\begin{array}{c}#1\\#2
\end{array}\right)}
\newcommand{\ga}{\alpha}
\newcommand{\gb}{\beta}
\newcommand{\gga}{\gamma}
\newcommand{\ul}{\underline}
\newcommand{\ol}{\overline}
\newcommand{\qed}{\kern 5pt\vrule height8pt width6.5pt depth2pt}
\newcommand{\Lrraro}{\Longrightarrow}
\newcommand{\Nb}{|\!\!/}
\newcommand{\NN}{{\rm I\!N}}
\newcommand{\bsl}{\backslash}
\newcommand{\gt}{\theta}
\newcommand{\op}{\oplus}
\newcommand{\C}{{\bf C}}
\newcommand{\Q}{{\bf Q}}
\newcommand{\Op}{\bigoplus}
\newcommand{\CR}{{\cal R}}
\newcommand{\tr}{\bigtriangleup}
\newcommand{\grr}{\omega_1}
\newcommand{\ben}{\begin{enumerate}}
\newcommand{\een}{\end{enumerate}}
\newcommand{\ndiv}{\not\mid}
\newcommand{\bab}{\bowtie}
\newcommand{\hal}{\leftharpoonup}
\newcommand{\har}{\rightharpoonup}
\newcommand{\ot}{\otimes}
\newcommand{\OT}{\bigotimes}
\newcommand{\bwe}{\bigwedge}
\newcommand{\gep}{\varepsilon}
\newcommand{\gs}{\sigma}
\newcommand{\rbraces}[1]{\left( #1 \right)}
\newcommand{\bbox}{$\;\;\rule{2mm}{2mm}$}
\newcommand{\sbraces}[1]{\left[ #1 \right]}
\newcommand{\bbraces}[1]{\left\{ #1 \right\}}
\newcommand{\OO}{_{(1)}}
\newcommand{\TT}{_{(2)}}
\newcommand{\FF}{_{(3)}}
\newcommand{\minus}{^{-1}}
\newcommand{\CV}{\cal V}
\newcommand{\CVs}{\cal{V}_s}
\newcommand{\un}{U_q(sl_n)'}
\newcommand{\on}{O_q(SL_n)'}
\newcommand{\slq}{U_q(sl_2)}
\newcommand{\olq}{O_q(SL_2)}
\newcommand{\UU}{U_{(N,\nu,\go)}}
\newcommand{\HH}{{\mathcal H}}
\newcommand{\ct}{\centerline}
\newcommand{\bs}{\bigskip}
\newcommand{\qua}{\rm quasitriangular}
\newcommand{\ms}{\medskip}
\newcommand{\noin}{\noindent}
\newcommand{\mat}[1]{$\;{#1}\;$}
\newcommand{\raro}{\rightarrow}
\newcommand{\map}[3]{{#1}\::\:{#2}\raro{#3}}
\newcommand{\alg}{{\rm Alg}}
\def\newtheorems{\newtheorem{theorem}{Theorem}[subsection]
                 \newtheorem{cor}[theorem]{Corollary}
                 \newtheorem{prop}[theorem]{Proposition}
                 \newtheorem{lemma}[theorem]{Lemma}
                 \newtheorem{defn}[theorem]{Definition}
                 \newtheorem{Theorem}{Theorem}[section]
                 \newtheorem{Corollary}[Theorem]{Corollary}
                 \newtheorem{Proposition}[Theorem]{Proposition}
                 \newtheorem{Lemma}[Theorem]{Lemma}
                 \newtheorem{Definition}[Theorem]{Definition}
                 \newtheorem{Example}[Theorem]{Example}
                 \newtheorem{Remark}[Theorem]{Remark}
                 \newtheorem{claim}[theorem]{Claim}
                 \newtheorem{sublemma}[theorem]{Sublemma}
                 \newtheorem{example}[theorem]{Example}
                 \newtheorem{remark}[theorem]{Remark}
                 \newtheorem{question}[theorem]{Question}
                 \newtheorem{Question}[Theorem]{Question}
                 \newtheorem{conjecture}{Conjecture}[subsection]}
\newtheorems
\newcommand{\proof}{\par\noindent{\bf Proof:}\quad}
\newcommand{\dmatr}[2]{\left(\begin{array}{c}{#1}\\
                            {#2}\end{array}\right)}
\newcommand{\doubcolb}[4]{\left(\begin{array}{cc}#1&#2\\
#3&#4\end{array}\right)}
\newcommand{\qmatrl}[4]{\left(\begin{array}{ll}{#1}&{#2}\\
                            {#3}&{#4}\end{array}\right)}
\newcommand{\qmatrc}[4]{\left(\begin{array}{cc}{#1}&{#2}\\
                            {#3}&{#4}\end{array}\right)}
\newcommand{\qmatrr}[4]{\left(\begin{array}{rr}{#1}&{#2}\\
                            {#3}&{#4}\end{array}\right)}
\newcommand{\smatr}[2]{\left(\begin{array}{c}{#1}\\
                            \vdots\\{#2}\end{array}\right)}

\newcommand{\ddet}[2]{\left[\begin{array}{c}{#1}\\
                           {#2}\end{array}\right]}
\newcommand{\qdetl}[4]{\left[\begin{array}{ll}{#1}&{#2}\\
                           {#3}&{#4}\end{array}\right]}
\newcommand{\qdetc}[4]{\left[\begin{array}{cc}{#1}&{#2}\\
                           {#3}&{#4}\end{array}\right]}
\newcommand{\qdetr}[4]{\left[\begin{array}{rr}{#1}&{#2}\\
                           {#3}&{#4}\end{array}\right]}

\newcommand{\qbracl}[4]{\left\{\begin{array}{ll}{#1}&{#2}\\
                           {#3}&{#4}\end{array}\right.}
\newcommand{\qbracr}[4]{\left.\begin{array}{ll}{#1}&{#2}\\
                           {#3}&{#4}\end{array}\right\}}

\title{{\bf Isocategorical Groups}} 
\author{Pavel Etingof\\
Massachusetts Institute of Technology\\
Department of Mathematics, Rm 2-165\\
Cambridge, MA 02139, USA, and\\
Columbia University\\
Department of Mathematics\\
2990 Broadway, New York, NY 10027, USA\\
{\rm email: etingof@math.mit.edu}
\and Shlomo Gelaki\\
Technion-Israel Institute of Technology\\
Department of Mathematics\\
Haifa 32000, Israel\\
{\rm email: gelaki@math.technion.ac.il}
}
\maketitle

\section{Introduction and the Main Results}
{\bf 1.} It is well known that 
the category $\text{Rep}(G)$ of 
finite dimensional complex representations of a finite group $G$, 
regarded as a symmetric tensor category, uniquely
determines the group $G$ up to isomorphism (see e.g. [DM, Theorem
3.2]).  

On the other hand, it is known that the Grothendieck ring of 
$\text{Rep}(G)$ (regarded as a ring with a distinguished basis, 
formed by the simple modules) is not sufficient to determine $G$.
For instance, the two nonisomorphic nonabelian groups of order
$8$ -- the group of symmetries of the square and the quaternion
group -- have the same Grothendieck rings. 

This raises the question
whether there is an intermediate amount of information 
between these two extremes that still determines $G$. 

Motivated by this question, we introduce the following notion 
of isocategorical finite groups. 

\begin{Definition} Two finite groups $G_1,G_2$ are called isocategorical if 
$\text{Rep}(G_1)$ is equivalent to $\text{Rep}(G_2)$ as a tensor category
(without regard for the symmetric structure).
\end{Definition}

The property of two groups to be isocategorical is much stronger 
than the property to have the same Grothendieck rings (with a
basis). For example, it is known (and will also follow from the
results below) that the two nonabelian groups of order 8 are not 
isocategorical. This raises the question whether isocategorical
groups must be isomorphic. 

Unfortunately (or fortunately, depending on the reader's taste), 
the answer to this question is negative. 
Namely, we prove the following result. 
\begin{Theorem}\label{th1} Let $Y$ be a vector space of dimension $\ge 3$ 
over a field of two elements. Let $V:=Y\oplus Y^*$ be the
corresponding symplectic vector space, and let
$\Sp(V)$ denote the
group of symplectic linear transformations of $V$. Then 
there exists a group $G_b$ which is isocategorical but not isomorphic
to $G:=\Sp(V)\ltimes V$. 
\end{Theorem}

The group $G_b$, constructed in Section 4, is in fact well known
in group theory. Namely, $G_b$ is an extension of 
$\Sp(V)$ by $V$ (using a cohomology class $b$), described by 
R. Griess in 1973 [Gr]; this extension can be defined 
over any field but is nontrivial only in characteristic $2.$ 

The group $G_b$ contains the {\em pseudosymplectic group} $\Ps(V)$ defined
by A. Weil in 1964 [W], and in particular the ``Weil
representation'' of $\Ps(V)$ (a projective representation on the
space of complex functions of $Y$) extends to $G_b$ (while 
by [LS] the group $G:=\Sp(V)\ltimes V$ and even $\Sp(V)$, for
$\dim(Y)\ge 4$, does not
admit nontrivial projective representations of such a small
dimension). 

We will call the group $G_b$ {\it the affine pseudosymplectic 
group}, since it is a ``modification'' of the {\em affine symplectic
group} $\ASp(V):=\Sp(V)\ltimes V$, and contains the pseudosymplectic
group
$\Ps(V)$. Correspondingly, we will denote $G_b$ by $\APs(V)$.   

{\bf 2.} In spite of Theorem \ref{th1},
 all groups isocategorical to a given group $G$ 
can be explicitly classified in group-theoretical terms. 
Their classification is described by the following theorem,
which is our second result. 

Let $A$ be a normal abelian subgroup of $G$
of order $2^{2m}$, and set $K:=G/A$. Let \linebreak
$R:A^\vee \to A$ be a $G$-invariant skew-symmetric isomorphism between $A$ and
its character group $A^\vee $ (i.e. $(Rx,x)=1\in \C^*$ for all $x\in A^\vee $).  
We can regard $R$ as a $G$-invariant nondegenerate 
skew-symmetric bilinear form 
on $A^\vee $ with values in $\C^*$. We will denote this form by the
same letter $R$: $(x,y)\mapsto R(x,y):=(Rx,y).$ 

The form $R$ defines a class in  
$H^2(A^\vee ,\C^*)^K$. Namely, this class is
represented by any 2-cocycle $J$ on $A^\vee $ with
values in $\C^*$, such that $R(x,y)=J(x,y)/J(y,x)$
(this class does not depend on the choice of $J$ since for any
two choices $J_1,J_2$ the cocycle $J_1J_2^{-1}$ is symmetric and
hence is a coboundary, see Lemma \ref{symcob} below).
We denote this class by $\bar R$.  

Let 
\begin{equation}\label{tau}
\tau:H^2(A^\vee ,\C^*)^K\to H^2(K,A)
\end{equation}
be the homomorphism defined as follows. 
For $c\in H^2(A^\vee ,\C^*)^K$, let $J$ be a $2-$cocycle representing
$c$. Then for any $g\in K$, the $2-$cocycle $J^gJ^{-1}$ 
is a coboundary (see 
Lemma \ref{symcob} below). Choose a cochain $z(g): A^\vee \to \C^*$ such that 
$dz(g)=J^gJ^{-1}$. Let 
\begin{equation}\label{btild}
\tilde b(g,h):=\frac{z(gh)}{z(g)z(h)^g}.
\end{equation}
It is easy to see that for any $g,h\in K$, 
the function $\tilde b(g,h):A^\vee \to \C^*$ is a character, i.e. 
$\tilde b(g,h)$ belongs to $A$
(see Lemma \ref{btildchar} below). Thus, $\tilde
b$ can be regarded
as a 2-cocycle of $K$ with coefficients in $A$
(our definition of the coboundary operator 
for group cohomology is the standard one, see [B], p.59). So
$\tilde b$ represents
a class $b$ in $H^2(K,A)$. It is easy to show that this class
depends only on $c$ and not on the choices we made. So we define
$\tau$ by $\tau (c)=b$. 

Let $b:=\tau (\bar R).$ Let $\tilde b$ be any cocycle representing $b$. 
For any $\gamma\in G$, let $\bar \gamma$ be the image of $\gamma$
in $K$. Introduce a new multiplication law $*$ on $G$ by
\begin{equation}\label{law}
\gamma_1*\gamma_2:=\tilde
b(\bar\gamma_1,\bar\gamma_2)\gamma_1\gamma_2.
\end{equation}
It is easy to show that this multiplication law introduces a 
new group structure on $G$, 
which (up to an isomorphism) depends only on $b$ and not on
$\tilde b$.  Let us call this group $G_b$.
\begin{Theorem}\label{th2} 
The following hold:
\ben
\item The group $G_b$ is isocategorical to $G$. 
\item Any group isocategorical to $G$ is obtained in this way. 
\een
\end{Theorem} 

Let us say that a finite group $G$ is {\em categorically rigid}
if any group isocategorical to $G$ is actually isomorphic to
$G$. 

\begin{Corollary}\label{th3} If $G$ is not categorically rigid then 
$G$ admits a normal abelian subgroup $A$, of order 
$2^{2m}$, $m\ge 1$, equipped with a skew-symmetric 
$G$-invariant isomorphism $R:A^\vee \to A$. In particular, all groups 
of orders $2k+1$, $2(2k+1)$, and all simple groups are
categorically rigid. 
\end{Corollary}
\begin{Remark} {\rm Note that Corollary \ref{th3} implies 
that the quaternion group is categorically rigid (all its normal
subgroups $A$ of order 4 are cyclic and hence do not admit 
a skew symmetric isomorphism $A^\vee \to A$). In particular, the two 
nonabelian groups of order 8 are not isocategorical. 
}
\end{Remark} 

The structure of the paper is as follows. 
In Section 2 we introduce some tools 
from Hopf algebra theory, which are needed to prove 
Theorems \ref{th1} and \ref{th2}. In Section 3, 
we prove Theorem \ref{th2}.
In Section 4, we prove Theorem \ref{th1}. 
In Section 5, we use the notion of isocategorical groups to
discuss the question when two triangular Hopf algebras are 
isomorphic as Hopf algebras. In Section 6, we discuss 
the properties of the affine pseudosymplectic group $\APs(V)$. 

We note that Theorem \ref{th2} and the results of section 5 
can be generalized without significant
changes to the case when finite groups are replaced by affine
proalgebraic groups. We do not discuss the details of this
generalization in this paper.

{\bf Acknowledgments.} We are indebted to R. Guralnick for many 
useful explanations and references on the group $\APs(V)$,
which were used in Section 6. We thank D. Kazhdan
for useful discussions. The first author is grateful
 to the IHES and the Erwin Schr\"odinger 
Institute for hospitality. The first author's work was conducted 
for the Clay Mathematical Institute. The second author is grateful 
to Bar Ilan University for its support.
The authors would also like 
to acknowledge the support of the NSF grant DMS-9700477. 

\section{Twists}
\subsection{Twists and Gauge Transformations}
Let $H$ be a  Hopf algebra. The multiplication, unit,
comultiplication, counit and antipode in $H$ will be
denoted by $m,1,\gD,\varepsilon,S$ respectively. We
recall Drinfeld's 
notion of a {\em twist} for $H$ [D]. 

\begin{Definition}\label{t}
A twist for $H$ is an invertible
element $J\in H\ot H$ which satisfies
\begin{equation}\label{t1}
(\Delta\ot I)(J)(J\ot 1)=(I\ot \Delta)(J)(1\ot J)\;\;
and \;\;
(\varepsilon\ot
I)(J)=(I\ot \varepsilon)(J)=1,
\end{equation}
where $I$ is the identity map of $H.$ 
\end{Definition}

Given a twist $J$ for $H,$ one can define a new
Hopf algebra
structure $(H^J,m,1,\Delta^J,\varepsilon,S^J)$ on the
algebra $(H,m,1)$ 
as follows. The coproduct is determined by
\begin{equation}\label{t2}
\Delta^J(a)=J^{-1}\Delta(a)J\;\text{for any}\;a\in H,
\end{equation}
and the antipode is determined by
\begin{equation}\label{tant}
S^J(a)=Q^{-1}S(a)Q\;\text{for any}\;a\in H,
\end{equation}
where $Q:=m\circ(S\ot I)(J).$ 

If $H$ is quasitriangular with the universal R-matrix $R$ then so 
is $H^J$
with the universal $R-$matrix $R^J:=J_{21}^{-1}RJ.$ 

If $J$ is a twist for $H$ and $x$ is an invertible
element of $H$ such that $\varepsilon(x)=1,$ then
\begin{equation}\label{jgauge}
J^x:=\Delta(x)J(x^{-1}\otimes x^{-1})
\end{equation}
is also a twist for $H.$
We will call the twists $J$ and $J^x$ {\em gauge
equivalent}, and $x$ a {\em gauge transformation}.
\begin{Remark} {\rm We note that in [EG2] the condition
$\varepsilon(x)=1$ was accidentally omitted and should be added.  
}
\end{Remark}

The map $(H^J,R^J)\raro (H^{J^x},R^{J^x})$
determined
by $a\mapsto xax^{-1}$ is
an isomorphism of quasitriangular Hopf algebras.
Thus, replacing a twist by a gauge equivalent one does 
not change the isomorphism class of $H^J$ as a quasitriangular
Hopf algebra. 
\subsection{Twists for Group Algebras of Abelian Groups}
Twists for group algebras of finite abelian groups will
be of particular interest in this paper. So let us give a well 
known description of such twists in terms of group $2-$cocycles. 

Let $A$ be a finite {\em abelian} group, and
$A^\vee :=\text{Hom}(A,\C^*)$ be its
character group. 

Then twists for
$H:=\C[A]$ are in one to one correspondence with 
$2-$cocycles $\tilde J$ of $A^\vee $ with coefficients in $\C^*,$
such that $\tilde J(0,0)=1$. 

Indeed, let $J$ be a twist for $H,$ and define
$\tilde J:A^\vee \times A^\vee \raro \C^*$ via
$\tilde J(\chi,\psi):=(\chi\ot\psi)(J).$ It is
straightforward to verify that $\tilde J$ is a $2-$cocycle
of $A^\vee $ (see e.g. [M, Proposition 3]), and that $\tilde
J(0,0)=1$.  

Conversely, let $\tilde J:A^\vee \times A^\vee \raro \C^*$ be a
$2-$cocycle of
$A^\vee $ with coefficients in $\C^*.$ For $\chi\in A^\vee ,$
let
$E_{\chi}:=|A|^{-1}\sum_{a\in A}\chi(a)a$ be the
associated idempotent of $H.$ Then
it is straightforward to verify that if $\tilde J(0,0)=1$ then 
$J:=\sum_{\chi,\psi\in A^\vee }\tilde J(\chi,\psi)E_{\chi}\ot
E_{\psi}$ is a twist for $H$ (see e.g. [M,
Proposition 3]). Moreover, it is easy to check that the above two assignments 
are inverse to each other.

 From now on we will abuse notation and denote twists for $\C[A]$ 
and corresponding $2-$cocycles by the same letter. 

 We note that gauge equivalence of twists for $H$
corresponds to
homological equivalence of the associated
$2-$cocycles. This fact will be useful later.

  In conclusion let us prove a lemma that is used many times
  throughout the paper. 

\begin{Lemma}\label{symcob}
A symmetric 2-cocycle $J$ of a finite abelian group $\Gamma$ 
with coefficients in $\C^*$ is a coboundary.
\end{Lemma}
\proof This is a simple special case of Corollary 3.3 of [EG2]. 
However, let us give an elementary direct proof.  

Recall that 2-cocycles encode abelian extensions.
A symmetric cocycle defines an extension that is itself an
abelian group.  
Therefore, the lemma claims that any exact sequence of abelian groups 
$$
0\to \C^*\to \tilde\Gamma\to \Gamma\to 0
$$
is split. Passing to the dual groups, we obtain the statement 
that any exact sequence of abelian groups
$$
0\to E\to \tilde E\to \Bbb Z\to 0,
$$
with $|E|<\infty$, is split. But this statement 
easily follows from the classification of finitely generated
abelian groups (i.e. from the fact
that $\tilde E$ decomposes into a direct sum of a
free part and the torsion part). \qed
\section{Proof of Theorem \ref{th2}}
Let us first prove part 2 of the theorem. 
We begin by reformulating the property of 
two groups to be isocategorical
in Hopf algebraic terms. This is accomplished by the following lemma. 

\begin{Lemma}\label{l1}
Finite groups $G_1$, $G_2$ 
are isocategorical if and only if there exists a twist $J$ for $\C[G_1]$ such that
$\C[G_1]^J$ and $\C[G_2]$ are isomorphic as Hopf
algebras (but {\em not}
necessarily as triangular Hopf algebras).
\end{Lemma}

\proof The "if" direction is clear, since 
twisting of a Hopf algebra does not change the tensor category 
of its representations (see e.g. [ES, Proposition 13.3]). To prove the 
"only if" direction, let $F:\Rep(G_1)\raro \Rep(G_2)$ be a tensor
equivalence, and let $F_2$ be the forgetful 
functor on $\Rep(G_2).$ Then $F_2\circ F$ is a
tensor (maybe nonsymmetric) 
fiber functor on $\Rep(G_1),$ so $\End(F_2\circ F)$
is the triangular Hopf algebra 
$\C[G_1]^J,$ where $J$ is some twist for $\C[G_1].$
On the other hand, $\End(F_2)=\C[G_2].$ 
So the functor $F$ establishes an isomorphism of
Hopf algebras 
$f:\C[G_1]^J\raro \C[G_2]$ (but not necessarily of
triangular Hopf algebras, since $F$ is not
assumed to be symmetric). We refer the
reader to [Ge] for more explanations. \qed

 So let us assume that 
$G_1$, $G_2$ are isocategorical,
and that such a twist $J$ has been fixed. 

\begin{Corollary}\label{cocomm}
The Hopf
algebra $\C[G_1]^J$ is cocommutative.
\end{Corollary}

\proof Clear.\qed

Let $H_J\subseteq \C[G_1]^J$ denote the span of 
(left or right) tensorands of the R-matrix
$R^J:=J_{21}^{-1}J$.
\begin{Proposition}\label{Rmat}
$H_J$ is a Hopf subalgebra of $\C[G_1]^J$, isomorphic to 
a group algebra of an abelian group. 
\end{Proposition}

\proof By the definition, $H_J$ is a minimal triangular 
Hopf algebra [R], i.e. 
it is generated by the components of its R-matrix. 
Therefore, $H_J$ is isomorphic, via the R-matrix, to its dual with 
opposite coproduct. But $H_J$ is cocommutative
by Corollary \ref{cocomm}, as it is a Hopf subalgebra
of $\C[G_1]^J$. Therefore, $H_J$ is commutative. We are done. \qed
 
\begin{Proposition} There exists a twist 
${\hat J}$ for $\C[G_1]$ such that $\C[G_1]^J$ is isomorphic 
to $\C[G_1]^{{\hat J}}$ as triangular Hopf algebras, and 
${\hat J}\in H_{{\hat J}}\otimes H_{{\hat J}}$. 
\end{Proposition}
\proof By construction, the Drinfeld element 
$u$ of the triangular Hopf algebra $(H_J,R^J)$ is $1$. 
Therefore, it is easy to show that there exists a twist $J'\in H_J\otimes H_J$ 
such that $R^J=R^{J'}$ (this is a very simple 
special case of [EG1, Theorem 2.1], which can be proved directly, without 
the use of Deligne's theorem).
Now, we see that 
$\C[G_1]^{J(J')^{-1}}$ is a triangular Hopf algebra
whose $R-$matrix is $1\otimes 1.$ Thus,
$J(J')^{-1}$ is a symmetric twist, hence it
is gauge equivalent to $1$ (see Corollary 3.3 of [EG2]). 
Thus, there exists an isomorphism 
of triangular Hopf algebras \linebreak $f:\C[G_1]^{J(J')^{-1}}\to 
\C[G_1]$. 
Therefore, $f:\C[G_1]^J\to \C[G_1]^{f(J')}$ 
is an isomorphism of triangular Hopf algebras.
It is clear that the twist ${\hat J}:=f(J')$ satisfies our 
requirements. \qed

Thus, we can assume, without loss of generality, that $J\in H_J\otimes H_J$. 
This implies that $H_J=\C[A]$, where $A$ is an abelian subgroup of $G_1$, and 
$J\in \C[A]\otimes \C[A]$. 
\begin{Proposition} $A$ is a normal subgroup of $G_1$, and 
the action of the group $K:=G_1/A$ on $A$ by conjugation preserves $R^J$.  
\end{Proposition} 
\proof By cocommutativity of $\C[G_1]^J$ we get 
$J^{-1}(g\ot g)J=J_{21}^{-1}(g\ot g)J_{21}$ for all
$g\in G_1.$
This implies that $R^J$ commutes with $g\ot g$ (here
we use that 
$A$ is abelian and hence $R^J=JJ_{21}^{-1}$). 
Since the left (and right) tensorands of $R^J$ span
$\C[A]$, the result follows. \qed

\begin{Lemma}\label{cor1} 
The twist $J$ can be chosen in such a way that $|A|=2^{2m},$ $m\ge 0.$
\end{Lemma}
\proof
Write $A=A_0\times A_1,$ where $A_0$ is the
$2-$Sylow
subgroup of $A,$ and the subgroup $A_1$ has odd order. 
Recall that twists for $\C[A]$ are $2-$cocycles of
$A^\vee $ with coefficients in $\C^*,$ and that 
gauge equivalence of twists corresponds to
homological equivalence of $2-$cocycles. It is well known that
$H^2(A,\C^*)=H^2(A_0,\C^*)\oplus
H^2(A_1,\C^*)$ (since the orders of $A_0,A_1$ are
relatively prime, see e.g. [CE]), 
so we can assume that $J=J_0J_1,$ where $J_0\in
\C[A_0]^{\ot 2},$ and $J_1\in\C[A_1]^{\ot 2}$. 
Consequently, the R-matrix $R:=R^J$ has a
factorization
$R=R_0R_1,$
where $R_0:=R^{J_0},$ and $R_1:=R^{J_1}.$ 

For all $x,y\in A_1,$ define $J_1'(x,y):=R_1(x/2,y)$
(this is well defined since $A_1$ has odd order). 
Let $R_1'$ be the $R-$matrix
associated with the twist $J_1'.$ Then
\begin{eqnarray*}
\lefteqn{
R_1'(x,y)=R_1(y/2,x)^{-1}R_1(x/2,y)=R_1(x,y/2)R_1(x/2,y)}\\
& = & R_1(x/2,y/2)^2
R_1(x/2,y)=R_1(x/2,y)^2=R_1(x,y).
\end{eqnarray*}
This implies that the $R-$matrix associated with
the twist $(J_1')^{-1}J_1$ is equal to $1\ot 1$. Hence, 
this twist is symmetric, and therefore 
is a coboundary (i.e. gauge equivalent to the identity),
by Lemma \ref{symcob}. 
We thus conclude that the twist
$(J_1')^{-1}J$ is gauge equivalent to $J_0.$ 

However, it is clear that $J_1'$ is $G_1-$invariant (since
$R_1$ is), hence $\C[G_1]^J$ and $\C[G_1]^{(J_1')^{-1}J}$
are isomorphic as Hopf
algebras, and the result follows. \qed

Thus, we will further assume that $|A|=2^{2m},$ $m\ge 0.$ 

Now we are ready to present the concluding part of the proof of 
part 2 of Theorem \ref{th2}. This part will
consist in calculating 
the group 
of grouplike elements of $\C[G_1]^J$, and showing that 
it is isomorphic to $(G_1)_b$ for an appropriate cohomology class $b$.

We will view $J$ 
not only as a twist but also 
as a $2-$cocycle of $A^\vee $ with values in $\C^*$, 
according to Section 2. 
For $g\in K$ let us write $J^g$ 
for the action of $g$ on $J$.  
Since $R^J$ is invariant under $G_1$, we get 
$$
J^gJ^{-1}=J_{21}^gJ_{21}^{-1}.
$$

This implies that the $2-$cocycle $J^gJ^{-1}$ of $A^\vee $ is symmetric, 
hence it is a coboundary, by Lemma \ref{symcob}.
Thus, there exists a function $z: K\to \C[A]^\times$, $g\mapsto z(g)$, 
(where $\C[A]^\times$ is the group of invertible elements of $\C[A]$),
such that 
$$
J^gJ^{-1}=\Delta(z(g))(z(g)^{-1}\otimes z(g)^{-1}).
$$ 

For $\gamma\in G_1$, let $\bar \gamma$ denote its image in $K$. 
Let us write $z(\gamma)$ for $z(\bar \gamma)$. We have 
$$
\Delta^J(\gamma)=J^{-1}J^{\bar \gamma}(\gamma\otimes \gamma)=
\Delta(z(\gamma))(z(\gamma)^{-1}\gamma\otimes 
z(\gamma)^{-1}\gamma). 
$$
This implies that $z(\gamma)^{-1}\gamma$ is a grouplike element in 
$\C[G_1]^J$. Thus, identifying $\C[G_1]^J$ with $\C[G_2]$, 
we obtain a map (of sets) $\varphi:G_1\to G_2$, 
$\varphi(\gamma)=z(\gamma)^{-1}\gamma$. 

It is easy to see that $\varphi$ is injective (and hence bijective). Indeed,  
$\varphi(\gamma_1)=\varphi(\gamma_2)$ if and only if
$z(\gamma_1)z(\gamma_2)^{-1}=\gamma_1\gamma_2^{-1}.$ Since 
$z(\gamma_1)z(\gamma_2)^{-1}\in
\C[A]$ and $\gamma_1\gamma_2^{-1}\in G_1$ we have that $\gamma_1=a\gamma_2$
for some $a\in A.$ But then, $z(\gamma_1)=z(a\gamma_2)=z(\gamma_2),$
hence $\gamma_1\gamma_2^{-1}=1.$ 

Finally, it is obvious from the definition of $\varphi$ that 
$$
\varphi(\gamma_1)\varphi(\gamma_2)=\tilde b(\bar\gamma_1,\bar\gamma_2)
\varphi(\gamma_1\gamma_2),
$$
where $\displaystyle{\tilde b(g,h)=z(gh)/z(g)z(h)^g\in \C[A]^\times.}$ 

\begin{Lemma}\label{btildchar}
The element $\tilde b:=\tilde b(g,h)$  is
a grouplike element, i.e. it belongs to $A$.
\end{Lemma}

\proof
$$
\Delta(\tilde b)(\tilde b^{-1}\otimes \tilde b^{-1})=
\frac{J^{gh}J^{-1}}{J^gJ^{-1}\cdot J^{gh}(J^g)^{-1}}=1.
$$
\qed

It is clear that $\tilde b$ is a $2-$cocycle of $K$ with
coefficients in $A.$
Let $b$ be the cohomology class of $\tilde b$ in $H^2(K,A)$.  
We have shown that 
$$
\varphi(\gamma_1)\varphi(\gamma_2)=\varphi(\gamma_1*\gamma_2),
$$
i.e. that $\varphi$ is an isomorphism $(G_1)_b\to G_2$. 
This completes the proof of part 2 of Theorem \ref{th2}, since
by the definition of $b$ we have $b=\tau (\ol {R^J}).$

Now let us prove part 1 of Theorem \ref{th2}. 
This part is essentially obvious from what we have done. Namely,
if $G$ is a finite group, 
$A$ its normal abelian subgroup, $K:=G/A$, and 
$b\in H^2(K,A)$ is given by  $b=\tau (\bar R),$ then choose a twist
$J\in \C[A]^{\otimes 2}$ such that $R=J_{21}^{-1}J$. Then, as we
have shown above, $\C[G]^J$ is isomorphic as a Hopf algebra to 
$\C[G_b]$, and hence by Lemma \ref{l1}, 
the groups $G$ and $G_b$ are isocategorical.  
Theorem \ref{th2} is proved. 
\begin{Remark} {\rm It is clear that in Theorem \ref{th2},
the ground field $\C$ (over which the tensor categories are
considered) can be replaced by any algebraically closed field 
of characteristic zero. On the other hand, one may introduce 
the following definition.}
\end{Remark}
\begin{Definition} Two finite groups $G_1,$ $G_2$ 
are $p-$isocategorical (for a prime $p$) 
if their tensor categories of representations
over an algebraically closed field of characteristic $p$ are
equivalent.  
\end{Definition}

Then we have the following.
\begin{Theorem}
Let $G_1,G_2$ be finite groups of order $n$ prime to $p.$
Then $G_1,G_2$ are isocategorical if and only if they are
$p-$isocategorical. In other words, the criterion for 
$G_1,G_2$ to be $p-$isocategorical is the same as in Theorem 
\ref{th2}.  
\end{Theorem}

The theorem is proved in the same way as Theorem \ref{th2}. 

\section{Proof of Theorem \ref{th1}}
Let $F$ be the field of two elements. 
Let $Y$ be an $n$-dimensional vector space over $F$, and 
$V:=Y\oplus Y^*$ be the corresponding symplectic space. 
Let $G:=\Sp(V)\ltimes V$. 
Let $A:=V$, and $K:=G/A=\Sp(V)$. 
Let $R:A^\vee \to A$ be the skew-symmetric isomorphism
defined by the symplectic form. 
Let $b:=\tau (\bar R)\in H^2(K,A)$.

The key step in the proof of Theorem \ref{th1} is the following 
result. 
\begin{Theorem}\label{nonzero} The element $b$ is nonzero if $n\ge 3$.
\end{Theorem}

Before proving Theorem \ref{nonzero}, let us 
show why it implies Theorem \ref{th1}. 
For this, it is enough to show that if $b\ne 0$ then $G_b$ is
not isomorphic to $G$ as a group. 

First of all, note that $K$ is a simple group (see [Go]).
This implies that $A$ is the unique normal subgroup of $G$ of 
order $|A|$ (indeed, if there is another one, it will project 
to a nontrivial normal subgroup of $K$, which is impossible). 
Thus, if $\psi: G\to G_b$ is an isomorphism, then $\psi(A)=A$. 
But this implies that $\psi$ is not only an isomorphism of 
groups but also an isomorphism of extensions
of $K$ by $A$, which means that 
$b$ must be zero, as desired.  

The rest of the section is devoted to the proof of 
Theorem \ref{nonzero}.

Let $P$ be the subgroup of $K$ which preserves 
the subspace $Y^*\subset V$. It is sufficient to show 
that the restriction of $b$ to $P$ is nonzero.
We will denote this restriction by $b_P$
(i.e. $b_P\in H^2(P,V)$). 
 
The group $P$ is naturally identified with 
the semidirect product $L\ltimes U$, where 
$L:=\GL(Y)$ and $U:=S^2Y^*$ is the additive group
of all self adjoint operators $Y\to Y^*$. 
The actions of $L$ and $U$ on $V$ are given by 
$l(v,f)=(lv,(l^*)^{-1}f)$, $l\in L$, and 
$u(v,f)=(v,f+uv)$, $u\in U$ (here $v\in Y$, $f\in Y^*$). 
   
Let us introduce a basis of $Y$: $e_1,...,e_n$. An element 
$v\in Y$ can be written as $v=\sum y_je_j$, where $y_j=y_j(v)$
are the coordinates of $v$. 

Define a function $\tilde\beta:P\times P\to Y^*$ by
the formula
\begin{equation}\label{beta}
\tilde\beta(u_1l_1,u_2l_2)=\chi(u_1,l_1u_2l_1^{-1}),
\end{equation}
where $\chi:U\times U\to Y^*$ is defined by the formula
\begin{equation}\label{chi}
\chi(u,u')(e_j)=(-1)^{\frac{1}{2}[(u+u')_{jj}-u_{jj}-u_{jj}']}, 
\end{equation}
in which $u_{ij}\in \lbrace{0,1\rbrace}$ are the matrix elements of
$u,$
and the terms on the right hand side 
(e.g. $u_{jj}$) are understood as {\bf integers}, 
0 or 1 (i.e. not as elements of $F$). 
\begin{Lemma} \label{l2} 
The following hold:
\ben
\item The element 
$\tilde\beta$ is a 2-cocycle of 
$P$ with coefficients in $Y^*$, representing some cohomology
class $\beta\in H^2(P,Y^*)$. 
\item The cohomology class $b_P$ is the image 
of $\beta$ under the 
map of cohomology induced by the embedding of 
$P$-modules $Y^*\hookrightarrow V$.
\een
\end{Lemma}
\proof Let us find a $2-$cocycle representing $b_P$. To do this, we pick $J$ 
to be
$J((v,f),(v',f'))=(-1)^{\lan v,f'\ran}$. 
Then $J^{ul}J^{-1}((v,f),(v',f'))=(-1)^{\lan v,uv'\ran}$ for all
$u\in U$ and $l\in L.$ 

Thus, we can take the function $z$ 
to be 
$$
z(ul)(v,f)=(-1)^{\sum_{m<j}u_{mj}y_my_j}\cdot i^{\sum_j u_{jj}y_j^2},
$$
where $i:=\sqrt{-1}$ (the exponent of $i$ is regarded as an
integer, 0 or 1). Therefore, 
we see that the function 
$$\tilde b(u_1l_1,u_2l_2)=
\frac{z(u_1l_1u_2l_2)}{z(u_1l_1)z(u_2l_2)^{u_1l_1}}$$
is given by
$$
\tilde b(u_1l_1,u_2l_2)=\chi(u_1,l_1u_2l_1^{-1}),
$$
where 
$$
\chi(u,u')(e_j)=i^{(u+u')_{jj}-u_{jj}-u_{jj}'},
$$
where the terms on the right hand side are understood as integers. 
This concludes the proof of the lemma. \qed
    
Now recall that we have a short exact sequence of $P$-modules
$$
0\to Y^*\to V\to Y\to 0
$$
(where $U$ acts trivially on both $Y$ and $Y^*$). 
To this sequence there corresponds a long exact sequence 
of cohomology. In particular, this long exact sequence 
has a portion 
\begin{equation}
H^1(P,Y)\to H^2(P,Y^*)\to H^2(P,V).
\end{equation}
Thus, by Lemma \ref{l2}, in order to show that $b_P\ne 0$ (and
hence to prove Theorem \ref{th1}), it is sufficient to prove the
following proposition. 

\begin{Proposition}\label{notimage} The element 
$\beta$ is {\bf not} contained in the image of $H^1(P,Y)$
under the connecting homomorphism
$H^1(P,Y)\to H^2(P,Y^*)$ of the long exact sequence. 
\end{Proposition}

The rest of this section is devoted to the proof of this proposition.

Let $\beta_U \in H^2(U,Y^*)$ be the restriction of $\beta$ to
$U$. 

\begin{Lemma} $\beta_U\ne 0$. 
\end{Lemma}

\proof It is easy to see from the explicit formula for 
$\tilde\beta,$ given in (\ref{beta}), that $\beta_U$ is
nonzero after restriction 
to the 1-dimensional subspace of $U$ generated by 
the element $u^{(11)}$, given by
$u^{(11)}_{mj}:=\delta_{1m}\delta_{1j}$ for all $m,j.$
\qed

Now we are ready to continue the proof of Proposition
\ref{notimage}.
Assume the contrary, i.e. that $\beta$ is the image of 
$\alpha\in H^1(P,Y)$. Let $\alpha_U\in H^1(U,Y)$ be the
restriction of $\alpha$ to $U$. Then $\beta_U$ is the image of $\alpha_U$
under the connecting homomorphism $H^1(U,Y)\to H^2(U,Y^*)$. 
Thus, to derive a contradiction, it is sufficient to show that 
$\alpha_U=0$.

It is clear that the element $\alpha_U$ is $L$-invariant, since it is
obtained by restricting a cohomology class of $P=L\ltimes U$
(this follows from the fact that the action of any group 
on its cochain complex 
with any coefficients descends to the trivial action on cohomology). 
Therefore, the statement that $\alpha_U=0$, and hence Proposition 
\ref{notimage}, follows from 

\begin{Lemma} If $n\ge 3$ then $H^1(U,Y)^L=0$. 
\end{Lemma}

\proof Since $U$ acts trivially on $Y$, we have
$H^1(U,Y)=\text{Hom}_F(U,Y)$. Since $F$ has two
elements, the space 
$U=S^2Y^*$ of symmetric bilinear forms on $Y$ is not an
irreducible representation of $L=\GL(Y)$. 
Namely, it has a submodule $\Lambda^2Y^*$
(skew-symmetric bilinear forms on $Y$, i.e. forms $B$ such that
$B(v,v)=0$ for all $v\in Y$), and the 
quotient is $Y^*$ (the map $S^2Y^*\to Y^*$ is given by 
$B\mapsto B',$ where $B'(v):=B(v,v)$ for all $v\in
Y$). It is easy to check 
(using explicit calculations with matrices) that the
representations $Y^*$, $\Lambda^2Y^*$ of 
$\GL(Y)$ are irreducible, and that 
$S^2Y^*$ is a nontrivial extension of one by the other. 
Therefore, an $L$-invariant homomorphism $U\to Y$ would have 
to be the pullback of an isomorphism $Y^*\to Y$. 
It is easy to show that such an isomorphism
of $L$-modules exists only for $n=1,2$. Thus, for $n\ge 3$ we have 
$H^1(U,Y)^L=0$, as desired. \qed

Proposition \ref{notimage}, Theorem \ref{nonzero}, and hence
Theorem \ref{th1} are proved. 

As we already mentioned in the introduction, 
we will call the group $G_b$ the {\em affine pseudosymplectic 
group}, and denote it by $\APs(V)$. 
The properties of this group are discussed in Section 6. 

\begin{Remark}{\rm It is easy to see that the construction 
of the group $\APs(V)$ makes sense over any finite field. 
However, it is easy to show that if the characteristic of the
ground field is odd then $\APs(V)$ is isomorphic to 
$\ASp(V):=\Sp(V)\ltimes V$. }
\end{Remark}

\begin{Remark} {\rm The cohomology class $\beta_U$, which appears 
in the proof of Theorem \ref{th1}, is invariant under $L=\GL(Y)$. 
However, our construction of this class used a basis of $Y$. 
So let us give a basis-free 
construction of this class.

Recall that if $U$ is a vector space over $F$, then $H^2(U,F)=
(S^2U)^*$. Thus, $H^2(U,Y^*)=(S^2U)^*\otimes Y^*=
\text{Hom}_F(S^2U,Y^*)$. Since we have an invariant projection 
$U\to Y^*$, we have an invariant projection 
$S^2U\to S^2Y^*$. Composing this projection with 
the invariant projection $S^2Y^*\to Y^*$, we get a
$\GL(Y)$ invariant map $S^2U\to Y^*$. One can show by a direct
calculation that this is exactly the cocycle $\beta_U$. 
}
\end{Remark}
\section{Hopf Algebra Isomorphisms of Triangular 
Semisimple Hopf Algebras}
Recall that in [EG2], we classified triangular semisimple 
Hopf algebras over $\C$, up to a triangular Hopf algebra 
isomorphism. 
In this section we will discuss a natural question which was not
discussed in [EG2]: when are
two triangular semisimple Hopf algebras over $\C$ isomorphic as usual
Hopf algebras (i.e. without regard for the triangular structure)? 
It turns out that the notion of isocategorical groups is useful 
in deciding this question. 

Let $B_1,B_2$ be two triangular semisimple Hopf algebras 
over $\C$. According to [EG1], 
$B_1=\C[G_1]^{J_1}$ and $B_2=\C[G_2]^{J_2}$ as 
Hopf algebras, where $G_i$ is a finite group, and $J_i$ 
is a twist for $\C[G_i],$ $i=1,2.$ Thus, in this section we will study 
the question: when are $\C[G_1]^{J_1}$,  
$\C[G_2]^{J_2}$ isomorphic as Hopf algebras? 

Recall that it follows from [EG1,EG2] that 
$\C[G_1]^{J_1}$,  
$\C[G_2]^{J_2}$ are isomorphic as {\bf triangular} Hopf algebras
if and only if there exists a group isomorphism 
$f:G_1\to G_2$ such that $f^{-1}(J_2)$ is gauge equivalent to $J_1$. 
To formulate a similar criterion for 
the rougher equivalence relation of isomorphism as usual Hopf
algebras, we need to generalize the notion 
of isomorphism of groups, and introduce the notion of
a categorical isomorphism. 
\begin{Definition} Let $G_1,G_2$ be finite groups. 
A categorical isomorphism $G_1\to G_2$ is a triple 
$(J,A,f)$, where 
$A$ is an abelian normal subgroup 
of $G_1$, $J\in \C[A]^{\otimes 2}$ 
a twist such that $R^J:=J_{21}^{-1}J$
is $G_1$-invariant and of maximal rank, 
and $f:\C[G_1]^J\to \C[G_2]$ an isomorphism 
of Hopf algebras.  
\end{Definition}

It is clear that two finite groups are isocategorical
if and only if there is a categorical isomorphism between them. 

There is an obvious notion of gauge equivalence 
of categorical isomorphisms, and for any $G_1,G_2$ 
there are finitely many categorical isomorphisms up to a gauge
equivalence. 

An example of a categorical isomorphism is 
a triple $(1,1,f)$, where $f$ is induced by an ordinary group 
isomorphism. Any categorical isomorphism which is
gauge equivalent to such will be called {\em trivial}. 
It is clear that for a trivial categorical isomorphism 
$(1,J,f)$, the map $f$ is actually an isomorphism of {\bf
triangular} Hopf algebras.
\begin{Theorem} Any isomorphism of Hopf algebras 
$\phi: \C[G_1]^{J_1}\to \C[G_2]^{J_2}$ 
is representable in the form $f\circ \psi$, 
where $(J,A,f):G_1\to G_2$ is a categorical isomorphism, and 
$\psi: \C[G_1]^{J_1}\to \C[G_1]^{Jf^{-1}(J_2)}$ is an 
isomorphism of triangular Hopf algebras. In particular, 
$\C[G_1]^{J_1}$,  
$\C[G_2]^{J_2}$ are isomorphic as Hopf algebras if and only if 
there exists a categorical isomorphism 
$(J,A,f):G_1\to G_2$, such that the twist
$Jf^{-1}(J_2)$ is gauge equivalent to $J_1$. 
\end{Theorem} 
\proof The map $\phi$ defines a Hopf algebra isomorphism
$\C[G_1]^{J_1\phi^{-1}(J_2)^{-1}}\to \C[G_2]$. 
By the results of Section 3, there exists a
categorical isomorphism $(J,A,f):G_1\to G_2$, 
and a triangular Hopf algebra isomorphism 
$\psi:\C[G_1]^{J_1\phi^{-1}(J_2)^{-1}}\to \C[G_1]^J$ such that 
$\phi=f\circ \psi$. The theorem is proved. \qed

Let us say that a finite group $G$ is {\it strongly categorically
rigid} if
any its normal abelian subgroup $A,$ which possesses 
a nondegenerate $G$-invariant 
skew-symmetric isomorphism \linebreak $R:A^\vee \to A$, is trivial. 
For example, any simple group is strongly categorically rigid.
\begin{Corollary} Suppose that $G_1$ or $G_2$ is strongly
categorically rigid. 
Then any isomorphism of Hopf algebras $\C[G_1]^{J_1}\to
\C[G_2]^{J_2}$
is actually an isomorphism of triangular Hopf algebras. 
\end{Corollary}
\proof
Assume that $G_1$ is strongly categorically rigid. 
Then any categorical isomorphism $(J,A,f):G_1\to G_2$ is trivial.
Therefore $f$ is a triangular Hopf algebra isomorphism, and the
result follows. \qed

\section{The Affine Pseudosymplectic Group}

The goal of this section is to discuss the properties of the
group $\APs(V)$, constructed in Section 4. These results are well
known and classical in the theory of finite groups; so the 
nature of this section is largely pedagogical. We note that much
of what is discussed below was explained to us by R. Guralnick. 

{\bf 1.} We start, following Weil [W],
by constructing a projective representation  
of $\ASp(V)$ called the ``Weil representation''. 
We will take the ground field $F$ to be the field of
$2$ elements, but the construction makes sense over any finite
field, yielding a projective representation of $\Sp(V)\ltimes V$ 
in odd characteristic. 

Let $\HH:=\text{Fun}(Y,\C)$ be the space of complex valued functions on $Y$. 
Consider the projective action $\rho$ of $V$ on $\HH$ given by 
$$
(\rho(y)f)(x):=f(x+y),\; (\rho(y^*)f)(x):=(-1)^{y^*(x)}f(x),\; y\in
Y,\;y^*\in
Y^*,\; \rho(y,y^*):=\rho(y^*)\rho(y). 
$$
Then 
$$
\rho(y_1,y_1^*)\rho(y_2,y_2^*)=\rho(y_1+y_2,
y_1^*+y_2^*)(-1)^{y_2^*(y_1)}.
$$
This shows that, 
$$
\rho(v_1)\rho(v_2)=\rho(v_1+v_2)J(v_1,v_2),
$$ 
where $J$ is the 2-cocycle defined in the proof of Lemma
\ref{l2}. In other words, $\rho$ extends to a representation 
of the ``Heisenberg'' group $E$, which consists of pairs 
$(v,c)$, $v\in V,c\in \lbrace \pm 1\rbrace$, with multiplication law
$(v_1,c_1)\cdot (v_2,c_2)=(v_1+v_2,c_1c_2J(v_1,v_2))$. 

Let $g\in \Sp(V)$. Consider the assignment $\rho^g$ defined by 
$\rho^g(v):=\rho(g^{-1}v).$ Then 
$$\rho^g(v_1)\rho^g(v_2)=\rho^g(v_1+v_2)J^g(v_1,v_2).$$

Let $z(g): V\to \C^*$ be as in the proof of Lemma \ref{l2}. 
Consider the map $\tilde \rho^g(v):=\rho^g(v)z(g)(v)$.
Then $\tilde \rho^g(v)$ also extends to a representation 
of $E$. 

It is easy to see that the representations 
$\rho,\tilde \rho^g$ of $E$ are isomorphic
(both are irreducible of dimension $2^{\dim(Y)}$, and such a
representation is unique). Therefore, there exists 
a unique, up to scaling, isomorphism 
$A(g):\HH\to \HH$, such that 
$A(g)\tilde\rho^g(v)=\rho(v)A(g)$. 

For a character $\chi: V\to \C^*$, let $T_\chi$ be an isomorphism 
$\HH\to \HH$, defined by 
$$
\rho(w)T_\chi=T_\chi\rho(w)\chi(w),\;w\in V.
$$
These operators exist and are unique up to scaling 
for the same reasons as $A(g)$. 

\begin{Proposition} One has 
$$
T_\chi T_\xi=c_1(\chi,\xi)T_{\chi+\xi},\;
A(g)T_\chi=c_2(g,\chi)T_{g\chi}A(g),
$$
$$
A(gh)=c_3(g,h)T_{\tilde b(g,h)}A(g)A(h),
$$
where $\tilde b$ is given in (\ref{btild}), 
and $c_1,c_2,c_3$ are suitable nonvanishing 
complex functions (depending on the choice of the normalization
of $A(g),T_\chi$). 
\end{Proposition}

\proof This follows from Schur's lemma, since the ratios of the
RHS and the LHS commute with $\rho(v)$. 
\qed

This shows that the assignment $A(\chi,g):=T_\chi A(g)$ 
defines a projective representation of the group $\APs(V)$ 
on $\HH$, which may be called the ``Weil representation'',
by analogy with [W]. 

\begin{Remark} {\rm The group $\Hom(V,\C^*)$ is naturally identified
with $V$ via $v\mapsto \chi_v$, $v\in V,$ where $\chi_v(w)=(-1)^{<v,w>}$. 
}
\end{Remark}

{\bf 2.} Here is another proof of Theorem \ref{th1} 
for $\dim(Y)\ge 4$. It was shown in [LS] that 
for $\dim(Y)\ge 4$, the group $\Sp(V)$
(and hence $\ASp(V)=\Sp(V)\ltimes V$) has no irreducible 
nontrivial projective representations of 
degree $2^{\dim(Y)}$ (or smaller). Thus, 
$\APs(V)$ is not isomorphic to $\ASp(V)$ (as it has the Weil
representation).

In particular, this argument shows that the statement 
that $\ASp(V)$ and $\APs(V)$ have equivalent representation categories,
does NOT extend to projective representations. 

{\bf 3.} Here is another construction of the group 
$\APs(V)$, given in [Gr]. 

Let $T$ be the group $(E\times \Z/4\Z)/\Gamma$, 
where $\Gamma$ is the subgroup of order $2$ generated by 
the element $(-1,2)$, 
$-1\in E$, $2\in \Z/4\Z$. Let $\text{Aut}_0(T)$ be the group of all
automorphisms of $T$ that act trivially on $\Z/4\Z$. 
Then $\text{Aut}_0(T)$ contains the subgroup $\text{Inn}(T)=T/Z(T)$ of inner
automorphisms of $T$, which is isomorphic to
$V$. It is easy to check (see   
[Gr]) that the group
$\text{Out}_0(T)=\text{Aut}_0(T)/\text{Inn}(T)$ is 
isomorphic to 
$\Sp(V)$. Thus, $\text{Aut}_0(T)$ is an extension of $\Sp(V)$ by $V$. 

The group $T$ has a unique irreducible representation $\HH$ 
of dimension $2^{\dim(Y)}$ with a fixed action of $Z(T)=\Z/4\Z$
(i.e. the generator acting by multiplication by $i$).
Therefore, the group $\text{Aut}_0(T)$ acts projectively in $\HH$. 

\begin{Proposition} The group $\text{Aut}_0(T)$ is isomorphic to 
$\APs(V)$, and the projective representation $\HH$ 
of $\text{Aut}_0(T)$ is the Weil representation.
\end{Proposition} 

\proof The proof is obtained directly from comparing the two
constructions of the Weil representation.
\qed

With this construction of $\APs(V)$, the essential part 
of the proof of Theorem \ref{th1} (i.e. the nontriviality of the 
extension) was proved in [Gr] (Corollary 2).

{\bf 4.} Recall [W] that the pseudosymplectic group 
$\Ps(V)$ is the group of all pairs $(g,Q)$,
where $Q$ is a quadratic form on $V$, and 
$g$ is a linear transformation of $V$, such that 
$$
Q(x+y)-Q(x)-Q(y)=J(gx,gy)-J(x,y), 
$$
with the operation $(g,Q)(g',Q')=(Q^{(g')^{-1}}+Q',gg')$. 
For fields of odd characteristic, this group is isomorphic to
$\Sp(V)$. 

Identify $T$ with $V\times \Z/4\Z$ in a natural way. 
The group $\Ps(V)$ acts on $T$ in the following way:
$$
(g,Q)(v,z)=(gv,z+2Q(v)).
$$
Thus, $\Ps(V)$ is a subgroup of $\APs(V)$. 
Namely, one can show that it is the preimage, under 
the projection $\APs(V)\to \Sp(V)$, of the orthogonal group 
$\text{O}(V)$ of the quadratic form $J(x,x)$. 
In particular, $\Ps(V)$ has the Weil representation,
constructed in [W] for local fields, and 
known even earlier for finite fields (see e.g. [BRW]).

\begin{Remark}{\rm In fact, as shown in [Gr],[BRW], 
the group $\Ps(V)$, for large enough $\dim(Y)$, is also 
a nontrivial extension (of $\text{O}(V)$ by $V$).}
\end{Remark}

\end{document}